\documentclass[11pt,leqno]{article} 
\usepackage{graphics}
\usepackage{xcolor}
\usepackage{lipsum}
\usepackage{subfigure}
\usepackage{graphicx}

\let\OLDthebibliography\thebibliography
\renewcommand\thebibliography[1]{
  \OLDthebibliography{#1}
  \setlength{\parskip}{2pt}
  \setlength{\itemsep}{2pt plus 0.3ex}
}

\newtheorem{thm}{Theorem}[section]
\newtheorem{lma}{Lemma}[section]

\newtheorem{prop}{Proposition}

\newcommand{\beqa}{\begin{eqnarray}}
\newcommand{\eeqa}{\end{eqnarray}}

\newcommand{\pf}{\noindent {\bf Proof:} $\s$ }
\newcommand{\epf}{ \hfill$\diamondsuit$ \medskip}

\newcommand{\beq}{\begin{equation}}
\newcommand{\eeq}{\end{equation}}
\newcommand{\lbl}{\label}
\newcommand{\s}{\; \;}

\newcommand{\om}{\omega}

\newcommand{\la}{\lambda}

\newcommand{\ra}{\rightarrow}
\newcommand{\al}{\alpha}

\title{Computation and analysis of global solution curves  for  super-critical equations}

\author{
Philip Korman  \thanks{Supported in part by the Taft faculty grant at the University of Cincinnati} \\ 
Department of Mathematical Sciences \\ 
University of Cincinnati \\ 
Cincinnati Ohio 45221-0025 \\
\\
Dieter S. Schmidt \\
Department of Computer Science \\
University of Cincinnati \\
Cincinnati, Ohio 45221-0030}

\date{}

\begin{document}

\maketitle
\begin{abstract} 
We study analytical and computational aspects for  Dirichlet problem on the unit ball $B$: $|x|<1$ in $R^n$, modeled on the equation
\[
\Delta u +\la \left(u^p+u^q \right)=0, \s \mbox{in $B$}, \s u=0 \s \mbox{on $\partial B$}, 
\]
with a positive parameter $\la$, and $1<p<\frac{n+2}{n-2}<q$, where $\frac{n+2}{n-2}$ is the critical power. It turns out that a special role is played by the Lin-Ni equation \cite{LN}, where $q=2p-1$ and $p>\frac{n}{n+2}$. This was already observed by I. Flores \cite{F}, who proved the existence of infinitely many ground state solutions. We study properties of  infinitely many solution curves of (\ref{i0}) that are separated by these ground state solutions. We also study singular solutions (where $u(0)=\infty$), and again the Lin-Ni equation plays a special role.
\medskip

Super-critical equations are very challenging computationally: solutions exist only for very large $\la$, and curves of positive solutions  make turns at very large values of $u(0)=||u||_{L^{\infty}}$. We overcome these difficulties by developing new results on singular solutions, and by using some delicate capabilities of {\em Mathematica} software. 
 \end{abstract}

\begin{flushleft}
Key words:  Computation of global solution curves, super-critical equations. 
\end{flushleft}

\begin{flushleft}
AMS subject classification: 65N25, 35J61.
\end{flushleft}

\section{Introduction}
\setcounter{equation}{0}
\setcounter{thm}{0}
\setcounter{lma}{0}

We study  positive solutions of the parameter depending Dirichlet problem on a unit ball $B$: $|x|<1$ in $R^n$, modeled on the equation
\beq
\lbl{i0}
\Delta u +\la \left(u^p+u^q \right)=0, \s \mbox{in $B$}, \s u=0 \s \mbox{on $\partial B$}, 
\eeq
with real $1<p<q$, and a positive parameter $\la$. Such equations arise in astrophysics, see the work of S. Chandrasekhar \cite{C}, for which he was awarded 1983 Nobel Prize in Physics.
By the classical theorem of B. Gidas, W.-M. Ni and L. Nirenberg positive  solutions of (\ref{i0}) are radially symmetric, so that $u=u(r)$, with $r=|x|$, and hence they satisfy 
\beq
\lbl{i1}
u''+\frac{n-1}{r}u'+\la \left(u^p+u^q \right)=0 \,,  \; u'(0)=0, \s  u(1)=0 \,.
\eeq
We consider both regular and singular solutions at $r=0$ (with $\lim _{r \ra 0^+} u(r)=\infty$ at singular solutions).
We are interested in super-critical equations with 
\beq
\lbl{i1+}
1 \leq p<\frac{n+2}{n-2}<q.
\eeq

By scaling we can study instead of (\ref{i1}) the positive solutions of  
\beq
\lbl{i2}
u''+\frac{n-1}{t}u'+u^p+u^q =0 \,, \; u'(0)=0, \; u(\omega)=0 \,,
\eeq
where $\om$ is the first root, and $\la=\om ^2$. If no finite $\om $ exists, we refer to the corresponding solution of (\ref{i2}) as a {\em ground state}.
Multiple solution curves of (\ref{i1}) are separated by ground state solutions. Indeed,  "shooting" for (\ref{i2}), starting with $u(0)>0$, results either with the first zero at some $\omega$ (solution of Dirichlet problem (\ref{i2})), or $u(r) \ra 0$ as $r \ra \infty$ (a ground state), while the case of a local minimum $u'(\xi)=0$ for some $\xi >0$ is clearly impossible.
\medskip

It turns out that super-critical problems have several distinguishing properties, compared with sub-critical problems that were mostly considered until recently.

\begin{enumerate}
  \item While in some cases there is exactly one solution curve (e.g., when $p=1$, as is shown below), the number of solution curves may be infinite (as in the case of Lin-Ni equation, as is seen below).
  \item A crucial role in determining the shape of the solution curve is played by the unique singular solution. While singular solutions exist for some sub-critical problems too (when $\frac{n}{n-2}<p<q< \frac{n+2}{n-2}$, see e.g., \cite{K4}), it appears that they do not play a big role in that case.
  \item Solution curves may have both horizontal and vertical asymptotes. In fact, the Lin-Ni equation has infinitely many horizontal asymptotes, and a vertical asymptote at some $\la _{\infty}>0$.
  \item Solution curves often make infinitely many oscillations around a vertical asymptote.
  \item Super-critical equations are very challenging computationally. It turns out that solution curves of (\ref{i0}) make turns at very large $\la$'s, while oscillations  around a vertical asymptote occur at very large values of the maximal value $u(0)$. For example, the solution curve for (\ref{i0}), with $n=3$, $p=4$, $q=7$, has a vertical asymptote at $\la _{\infty} \approx e^{11.44} \approx 92967$.
\end{enumerate}

It turns out that a special role is played by the  Lin-Ni equation
\beq
\lbl{-i2}
u''+\frac{n-1}{r}u'+u^p+u^{2p-1}=0 \,, \s u'(0)=0 \,, \; u(1)=0 \,,
\eeq
with $\frac{n}{n-2}<p<\frac{n+2}{n-2}$, which implies that $2p-1>\frac{n+2}{n-2}$. It was discovered in Lin-Ni \cite{LN} that this equation has an explicit ground state solution
\[
u(r)=\left(\frac{2p(np-n-2p)}{(np-n-2p)^2+p(p-1)^2r^2}\right)^{\frac {1}{p-1}}. 
\]
Then I. Flores \cite{F} discovered that the  Lin-Ni equation has infinitely many other ground state solutions, and it is the only equation in the class (\ref{i0}) with this property.
 It turns out  that there is  a discrete infinite sequence of ground states, which serve to separate an infinite sequence of solutions to the Dirichlet problem. These ground states are of fast decay, and they tend to the Lin-Ni ground state, which is of slow decay. We present computational results on the global bifurcation diagram for the Dirichlet problem, which use some delicate {\em Mathematica} capabilities.  Computing infinitely many solution curves in Figure \ref{fig:3} below involved a tremendous challenge: one needed to go to astronomically high $\la$'s, and then zoom in on a small neighborhood of $u(0)=2$.
\medskip

Our results also suggest that the Lin-Ni equation is very special. We show that the  Lin-Ni equation has a singular solution of the form $u(r)=r^{-\frac{1}{p-1}}v(r)$ where $v(r)$ is an analytic function. For all other equations of type (\ref{i0}) the results on singular solutions are similar, but $v(r)$ is not analytic at $r=0$. We also mention that some kind of reduction of order occurs for the  Lin-Ni equation. It was shown in P. Korman \cite{K2} that the Lin-Ni ground state solution of (\ref{-i2}) satisfies the first order equation
\beq
\lbl{-i20}
u'(r)=-Aru^{p},
\eeq
where $A=\frac{p-1}{np-n-2p}$. (Observe that $A>0$ iff $n>\frac{n}{n-2}$, and then $2p-1>\frac{n+2}{n-2}$, the critical exponent.) The Lin-Ni solution (of both (\ref{-i2}) and (\ref{-i20})) can be written as $u(r)=\left(\frac{An-1}{1+pA^2r^2}\right)^{\frac {1}{p-1}}$.
\medskip

\begin{figure}[ht!]

\subfigure[]{\includegraphics[scale=0.6]{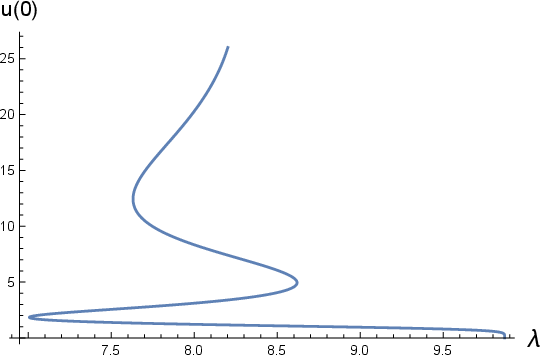}}
\hfill
\subfigure[]{\includegraphics[scale=0.75]{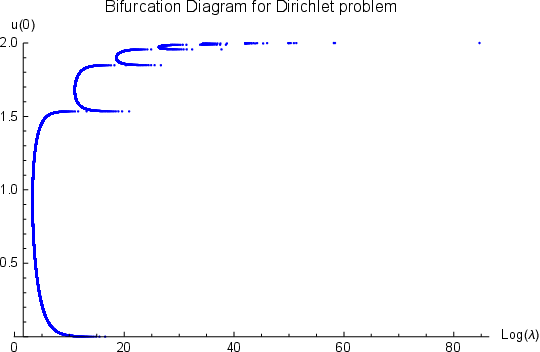}}
\caption{ a)\ The unique solution curve of (\ref{i1}), $n=3$, $p=1$, $q=7$. \newline  b)\ Infinitely many solution curves  of (\ref{i1}), $n=3$, $p=4$, $q=7$.}
\label{fig:1a}
\end{figure}

Our computations show that both the smaller and larger powers in (\ref{i0}) play an important role in determining the global bifurcation diagram. In Figure \ref{fig:1a} we see how the change in the smaller power $p$ produces a profound change of the bifurcation diagram.

\medskip

Both analytical and computational results of this paper extend to $p$-Laplacian equations, by proceeding as in \cite{KS1} and \cite{K18}, and to non-autonomous equation of the type
\[
\Delta u +\la r^{\al} \left(u^p+u^q \right)=0, \s \mbox{in $B$}, \s u=0 \s \mbox{on $\partial B$},
\] 
by using the change of variables in \cite{K2} that removes the $r^{\al}$ term while preserving the form of the equation and boundary conditions (while changing the dimension $n$).

\section{Uniqueness and non-uniqueness of solution curves}
\setcounter{equation}{0}
\setcounter{thm}{0}
\setcounter{lma}{0}

In case the nonlinear term $f(u)$ changes sign for $u>0$, multiple  curves of positive solutions are typical, see e.g., P. Korman \cite{K1}.
If $f(u)>0$ for $u>0$, multiple solution curves are also possible, although only for supercritical problems. It appears unclear for which  supercritical $f(u)$ multiple solution curves occur. We begin with a case when solution curve is unique.

\begin{thm}\lbl{prop:1}
If $p=1$ and $q>1$ is arbitrary, all positive solutions of (\ref{i1}) lie on a unique solution curve in the $(\la ,u(0))$-plane.
\end{thm}

\pf
Let us recall the following well known facts, see e.g.,   P. Korman \cite{K1}. There is a continuous curve of positive solutions bifurcating from zero at $\la =\la _1$. If this curve extends to $u(0)=\infty$, this solution curve is unique, since $u(0)$ is a global parameter uniquely identifying the solution pair $(\la, u(r))$. Indeed, all possible values of $u(0)$, $0<u(0)<\infty$ are taken up by this curve.
\medskip

If $u(0)$ remains bounded on this curve, then the increasing values of $u(0)$ tend to a finite limit, call it $\bar u$. Along this curve $\la \ra \infty$, since it is known that all curves of positive solutions extend globally, see e.g., P. Korman \cite{K1}. Let $v(t)$ be the solution of (\ref{i2}), with $v(0)=\bar u$. We claim that  $v(t)$ is a ground state, i.e., a decreasing solution tending to zero as $t \ra \infty$. Indeed, near $v(t)$ we have solutions of (\ref{i1}) with arbitrarily large $\la$, which correspond after rescaling to solutions of
\beq
\lbl{1}
u''+\frac{n-1}{t}u'+u+u^q =0,   \; u'(0)=0, \; u(\omega)=0
\eeq
with arbitrarily large $\omega=\sqrt{\la}$. If $v(t)$ were to vanish at some finite $t$, we would have a contradiction with the continuous dependence of solutions on initial data. We show next that the ground state solution $v(t)$ is not possible.
\medskip

Near the point of bifurcation from zero (which occurs for (\ref{i2}) at $\la=\la _1$, the principal eigenvalue of the Laplacian on the unit ball, with zero boundary condition) we have small solutions $u(r)$ of (\ref{i1}), which correspond to small solutions $u(t)$ of (\ref{i2}), with $\om \approx \sqrt{\la _1}$:
Select one such solution, with 
\beq
\lbl{2}
u(t)<v(t) \s \mbox{ on $(0,\om )$}. 
\eeq
So that $v(t)$ satisfies
\[
v''+\frac{n-1}{t}v'+v+v^q =0,   \; v'(0)=0, \; v(\omega)>0,
\]
or
\beq
\lbl{3}
\left(t^{n-1}v'  \right)'+t^{n-1}(v+v^q )=0,   \; v'(0)=0, \; v(\omega)>0.  
\eeq
Multiply (\ref{3}) by $u$, and subtract the equation (\ref{1}) multiplied by $v$, then integrate over $(0,\om)$. Obtain
\[
-\om ^{n-1} u'(\om)v(\om)+\int_0^{\om} t^{n-1} uv \left(v^{q-1}-u^{q-1} \right) \, dt=0.
\]
In view of (\ref{2}), both terms on the left are positive, a contradiction.
\epf

\noindent
{\bf Remark } Our computations indicate that the solution curve is still unique for $p>1$, and $p \approx 1$. It is an interesting open question to determine a critical $p>1$ at which the solution curve breaks into pieces. The argument above does not work for $p>1$ since small solutions of (\ref{i2}) may have large first roots. 
\medskip

\noindent
{\bf Example }
When $n=3$ and $q=7$ the equation (\ref{i1}) takes the form
\[
u''+\frac{2}{r}u'+\la \left(u^p+u^7 \right)=0 \,,  \; u'(0)=0, \s  u(1)=0 \,.
\]
Our numerical computations show that the solution curve is unique for $1 \leq p \leq 3.7$, and it breaks into two pieces by $p=3.8$ (here $p=4$ corresponds to the Lin-Ni equation). As $p$ is increased toward  $4$, the number of solution curves increases, and it is infinite at $p=4$, see Figure \ref{fig:1a} (b).

\begin{prop}\lbl{prop:2}
In case  $1<p<\frac{n+2}{n-2}<q$, there is a number $A>0$ so that all positive solutions of (\ref{i1}) satisfying $u(0)>A$ lie on a unique continuous solution curve in the $(\la ,u(0))$-plane.
\end{prop}

\pf
Multiple curves of positive solutions of (\ref{i1}) are separated by ground states. We will show later on that  ground states with large $u(0)$ are not possible. (Solutions of (\ref{i2})  with large $u(0)$ vanish at $\om \approx \sqrt{\la ^*}$, with $\la ^*$  corresponding to the singular solution of (\ref{i1}).)
\epf

Let us also recall the following known result, see K. Schmitt \cite{S20} or P. Korman  \cite{K1}. It can also be derived from the above propositions, and the fact that solution curves approach a singular solutions at some $\la ^*$ as $u(0) \ra \infty$.

\begin{prop}\lbl{prop:2a}
In case  $1 \leq p<\frac{n+2}{n-2}<q$, there are numbers $0<\bar \la <\bar \Lambda$ so that $\bar \la< \la < \bar \Lambda $ holds for all positive solutions of (\ref{i1}).
\end{prop}

For subcritical problems the situation is much simpler, see e.g., P. Korman \cite{K1} or T. Ouyang and J. Shi \cite{OS}.

\begin{prop}
In case $1<p,q<\frac{n+2}{n-2}$  the problem (\ref{i1}) has a unique positive solution for each $\la >0$. These solutions lie on a  unique solution curve, which is hyperbola-like, with $u(0) \ra \infty$ as $\la \ra 0$, and
$u(0) \ra 0$ as $\la \ra \infty$. 
\end{prop}

\section{Singular solution for the power case}
\setcounter{equation}{0}
\setcounter{thm}{0}
\setcounter{lma}{0}

In this section we collect the relevant results on the Lane-Emden equation (which arose in astrophysics in the 19-th century)
\beq
\lbl{100}
w''+\frac{n-1}{r}w'+w^p=0
\eeq
that are more or less known, and apply them to the solution curves of (\ref{i0}). The equation possesses (for $p>\frac{n}{n-2}$) a   {\em singular solution}, which is unbounded at $r=0$. Singular solution plays an important role in determining the shape of solution curves for the Dirichlet problems (\ref{i0}). We describe the singular solution next.
\medskip

Define 
\beq
\lbl{15}
B=\frac{2}{(p-1)^2} \left[-2+(p-1)(n-2)   \right]. 
\eeq
Observe that $B>0$ if and only if $p>\frac{n}{n-2}$. The following lemma is proved by a direct calculation.

\begin{lma}\lbl{lma:lane}
Assume that $p>\frac{n}{n-2}$. Then $w_0(r)=B^{\frac{1}{p-1}}r^{-\frac{2}{p-1}}$ is a singular ground state solution of (\ref{100}).
\end{lma}

Linearization of the equation (\ref{100}) around the singular solution $w_0(r)$ is
\[
z''+\frac{n-1}{r}z'+pw_0^{p-1}z=0 \,, \s z'(0)=z(1)=0 \,,
\]
which is the same as 
\beq
\lbl{6}
z''+\frac{n-1}{r}z'+\frac{pB}{r^2}z=0 \,, \s z'(0)=z(1)=0 \,.
\eeq
This is  Euler's equation. Its characteristic equation $s^2+(n-2)s+pB=0$ has the roots
\[
s=\frac{-(n-2) \pm \sqrt{(n-2)^2-4pB}}{2}\,.
\]
The roots are complex, so that $z(t)$ is oscillatory, provided that
\beq
\lbl{7}
(n-2)^2-4pB<0 \,.
\eeq

To analyze the inequality (\ref{7}), we begin with the equal sign:
\beq
\lbl{17}
(n-2)^2-4p\frac{2}{(p-1)^2} \left[-2+(p-1)(n-2)   \right]=0.
\eeq
Solving this quadratic equation for $p$ gives the roots $p_1,p_2$, or two solution curves $p_1(n),p_2(n)$
\beq
\lbl{15a}
p_1(n)=\frac{n^2 - 8 n-8
   \sqrt{n-1}+4}{n^2-12 n+20},
\eeq
\beq
\lbl{15b}
p_2(n)=\frac{n^2 - 8 n+8
   \sqrt{n-1}+4}{n^2-12 n+20}.
\eeq

By an elementary analysis one shows that, see Figure \ref{fig:0a},
\beqa
\lbl{18}
&  p_2(n)<0<p_1(n)<\frac{n+2}{n-2}, \s\s \mbox{for $2<n<10$} \\
&  1<p_1(n)<\frac{n+2}{n-2}<p_2(n), \s\s \mbox{for $n>10$}. \nonumber
\eeqa

\begin{lma}\lbl{lma:4}
If $p>\frac{n+2}{n-2}$ and $2<n \leq 10$, then the inequality (\ref{7}) holds.
\end{lma}

\pf
If $2<n<10$, the numerator in (\ref{17}) is a quadratic in $p$ with the leading coefficient (in $p^2$)
\beq
\lbl{19}
(n-2)^2-8(n-2)=(n-2)(n-10)<0,
\eeq
and hence this  quadratic  is negative when $p$ is greater  than its larger root. The root $p_2(n)$ is negative, while the larger root satisfies $0<p_1(n)<\frac{n+2}{n-2}$, so that  (\ref{7}) holds, since $p>\frac{n+2}{n-2}>p_1(n)$.
\medskip

When $n=10$, $\frac{n+2}{n-2}=\frac{3}{2}$,
\[
(n-2)^2-4 p B=\frac{64-48 p}{(p-1)^2}<0,
\]
for $p>\frac{3}{2}$.
\epf

\begin{lma}\lbl{lma:5}
If $p>\frac{n+2}{n-2}$ and $n> 10$, then the inequality (\ref{7}) holds if and only if $p<p_2(n)$.
\end{lma}

\pf
The numerator in (\ref{17}) is a quadratic in $p$, and the leading coefficient is now positive in view of (\ref{19}). This quadratic is negative when $p$ is between its roots. 
\epf

\medskip

\begin{figure}
\begin{center}
\scalebox{0.9}{\includegraphics{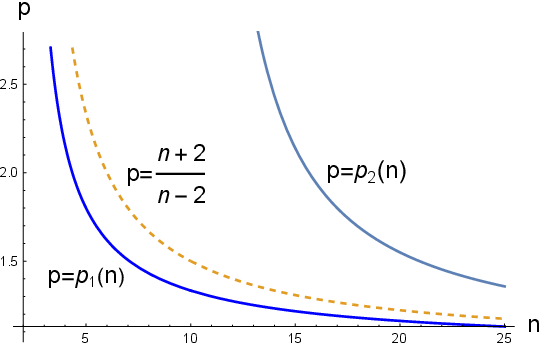}}
\end{center}
\caption{ The  critical curves }
\lbl{fig:0a}
\end{figure}

The number $p_2(n)$ is known as the Joseph-Lundgren exponent, introduced in \cite{JL}. It can be written as 
\beq
\lbl{JL}
p_2(n)=1+\frac{4}{n-2 \sqrt{n-1}-4}.
\eeq
Define $ p_{_{JL}}=\infty$ for $3 \leq n \leq 10$ and  $ p_{_{JL}}=1+\frac{4}{n-2 \sqrt{n-1}-4}$ for $n \geq 11$. The following theorem  was included in  Y. Miyamoto \cite{M}.

\begin{thm}(\cite{M})
Assume that $1<p< \frac{n+2}{n-2} <q< p_{_{JL}}$. Assume that the singular solution of  (\ref{i0}) occurs at $\la=\la _{\infty}$. Then the solution curve of (\ref{i0}) makes infinitely many turns around $\la=\la _{\infty}$.
\end{thm}

Our computations indicate that  the solution curve actually has a vertical asymptote at  $\la=\la _{\infty}$, and that the value of $\la _{\infty}$ is affected drastically by the lower power of $p$. 
\medskip

Let us recall Theorem 4.1 in M.-F. Bidaut-Veron \cite{B}.
\begin{thm}\lbl{lma:4*} (\cite{B})
Any solution of the equation  (\ref{100}) is either of slow decay (when $w(r) \sim cr^{-\frac{2}{p-1}}$ as $r \ra \infty$), or of fast decay as $r \ra \infty$ (when $w(r) \sim cr^{-\frac{1}{n-2}}$ as $r \ra \infty$).
\end{thm}

Let now $w(r)$ be the solution of 
\beq
\lbl{60}
w''+\frac{n-1}{r}w'+w^p=0 \,, \s w'(0)=0 \,, \; w(0)=1.
\eeq
It is well-known (by Pohozhaev's identity) that  for super-critical $p>\frac{n+2}{n+2}$, $w(r)$ remains positive for all $r \in (0,\infty)$, it is decreasing for all $r$, and tends to zero as $r \ra \infty$, i.e., $w(r)$ is a {\em ground state solution}.

\begin{lma}\lbl{lma:slow}
The ground state solution $w(r)$ of (\ref{60}) is of slow decay.
\end{lma}

\pf
By the Theorem 5.1 part (iv) in  M.-F. Bidaut-Veron \cite{B} the possibility of  fast decay is ruled out by the fact that the solution does not oscillate at $r=0$. (I. Flores \cite{F} has an alternative proof.)
\epf

The following lemma of S. Chandrasekhar \cite{C}, which gives the exact form of the slow decay, was used in D.D. Joseph and T.S. Lundgren \cite{JL} and later in C. Budd and J.  Norbury \cite{BN}. We  present a short proof, based on work of M.-F. Bidaut-Veron \cite{B}, used in the preceding lemma. We shall use this lemma later to calculate the value of $\la _{\infty}$ for the problem (\ref{+*}) below.

\begin{lma}\lbl{lma:50}
For $p>\frac{n+2}{n-2}$, the solution  of (\ref{60}) satisfies $\lim_{r \ra \infty} r^{\frac{2}{p-1}}w(r)=B^{\frac{1}{p-1}}$ ($B$ was defined in (\ref{15})).
\end{lma}

\pf
By Lemma \ref{lma:slow}, $w(r)$ is of slow decay as $r \ra \infty$. Make a change of variables $s=\frac{1}{r}$, and study the solution near $s=0$. Obtain
\beq
\lbl{n1}
s^4w''+(3-n)s^3w' +w^p=0 \,.
\eeq
Make a change of variables 
\beq
\lbl{n2}
w(s)=s^{\frac{2}{p-1}} v(s) \,.
\eeq
Upon dividing by $s^{\frac{4 p+2}{p-1}}$, the equation (\ref{n1}) takes the form
\beq
\lbl{n3}
v''+\frac{c_1}{s}v'+\frac{1}{s^2}g(v)=0 \,,
\eeq
where $c_1=\frac{-n p+n+3 p+1}{p-1}$, and $g(v)=\frac{4 p-2 n (p-1)}{(p-1)^2}v+v^p$. By Lemma \ref{lma:slow}, $v(0)$ is finite. From the equation (\ref{n3}), we see that 
\beq
\lbl{n4}
g(v(0))=0 \,,
\eeq
since otherwise $v(s)$ tends to infinity at $s=0$, which follows by integration of (\ref{n3}), and contradicts Lemma \ref{lma:slow}. Solving (\ref{n4}) one gets $v(0) =B^{\frac{1}{p-1}}$, completing the proof, since $v(0)=\lim_{s \ra 0} s^{-\frac{2}{p-1}}w(s)=\lim_{r \ra \infty} r^{\frac{2}{p-1}}w(r)$.
\epf

\section{The global picture for the Lin-Ni equation}
\setcounter{equation}{0}
\setcounter{thm}{0}
\setcounter{lma}{0}

In this section we present surprising computational results on the global solution curves for Dirichlet problem for the Lin-Ni equation
\beq
\lbl{d1}
u''+\frac{n-1}{r}u'+\la \left(u^p+u^{2p-1} \right)=0 \,, \s u'(0)=0 \,, \; u(1)=0 \,,
\eeq
with $\frac{n}{n-2}<p<\frac{n+2}{n-2}$, which implies that $2p-1>\frac{n+2}{n-2}$, and $\la $ is a positive parameter. At the end of this section we provide an analytical support to our computational results, by studying the linearized problem for (\ref{d1}).

\medskip

A simple scaling 
argument shows that the value of $u(0) > 0$ is a global parameter, uniquely
identifying the solution pair $(\la , u(r))$ of (\ref{d1}). Therefore we shall study the solution curves $(\la , u(0))$. Points on this curve are numerically computed as follows: given any $u(0)=u_0>0$, solve the initial value problem 
\beq
\lbl{d2}
u''+\frac{n-1}{r}u'+ u^p+u^{2p-1} =0 \,, \s  u(0)=u_0 \s u'(0)=0  \,,
\eeq
and let $r>0$ be the first root of $u(r)$. Then $\la =r^2$, as easily follows by scaling, or see e.g., P. Korman \cite{K1}.
\medskip

    We present our computations in case $n=3$ and $p=4$. The Mathematica program below is one way of computing the solution curves $(\la , u(0))$ of (\ref{d1}).  We also have done the computations with higher  precision. In most cases the standard double floating point arithmetic seems to be sufficient. Next we present a sample {\em Mathematica} program (with detailed explanations given in our notebook published  in Wolfram Notebook Archive \cite{KS}). The crucial element involves using the WhenEvent option of the NDSolve command to stop integration when the first root of solution is reached.
    
    \begin{verbatim}
rep={n->3,p->4}
f[u_]=u^p+u^(2p-1)
dequ=u''[r]+(n-1)/r*u'[r]+f[u[r]]/.rep
umin=0;
umax=10000000;
nsteps=2000;
delu=(umax-umin)/nsteps;
h=10^-100;
dirichlet={};
For[k=1,k <= nsteps,k++,
  u0= umin+k*delu ;
  b0=True;
  sol1=NDSolveValue[{dequ==0,u[h]==u0,u'[h]==0,
  WhenEvent[u[r]== 0,
  b0=False;
  AppendTo[dirichlet,{r^2,u0}];
  "StopIntegration"]},
   u,{r,h,1/h}];
If[b0,
  Print["no zero for u0=",u0];
];

ListPlot[dirichlet,Joined -> False,  AxesLabel -> {"\[Lambda]","u(0)"}, 
PlotRange -> {{280, 325},{umin,umax}}, PlotStyle ->{Blue, Thickness[0.01]},
PlotLabel->"Bifurcation Diagram for Dirichlet problem"]
    \end{verbatim}
\medskip

We now describe the global picture in terms of the global parameter $u(0)$. For $0<u(0)<2$ there is an increasing sequence $\{u_i \} \ra 2$, and an infinite sequence of parabola-like solution curves in the $(\la , u(0))$ plane that are opening to the right. The first of these curves has a unique turning point at $\la =\la _1 \approx 735$. As $\la \ra \infty$, $u(0) \ra 0$ along the lower branch, while $u(0) \ra u_1 \approx 1.54$ along the upper branch. The curve $i$, $i \geq 2$,  has a unique turning point at $\la =\la _i$. (Observe that $\la _2 \approx 3.5 \times 10^9>>\la _1$, and $\la _2$ is a very large number.) As $\la \ra \infty$, $u(0) \ra u_{i-1}$ along its lower branch, while $u(0) \ra u_i$ along the upper branch. The  sequence $\{\la _i \}$ is increasing (very rapidly), while the  sequence $\{u _i -u_{i-1}\}$ is decreasing, see Figure \ref{fig:1a} (b).
\medskip

Turning to the $u(0)>2$ range, there is a solution curve along which $u(0) \downarrow v_1  $ as $\la \ra \infty$. This curve makes infinitely many turns, and $u(0) \ra \infty$, as $\la $ tends to a finite $\la _{\infty}$. For $2<u(0)<v_1$ there is an decreasing sequence $\{v_i \} \ra 2$, and an infinite sequence of parabola-like solution curves in the $(\la , u(0))$ plane that are opening to the right. The curve $i$, $i \geq 1$,  has a unique turning point at $\la =\mu _i$. As $\la \ra \infty$, $u(0) \ra v_{i+1}$ along its lower branch, while $u(0) \ra v_i$ along the upper branch. The  sequence $\{\mu _i \}$ is increasing, while the  sequence $\{v _i -v_{i+1}\}$ is decreasing. 
\medskip

\begin{figure}[htb!]
\scalebox{1.2}{\includegraphics{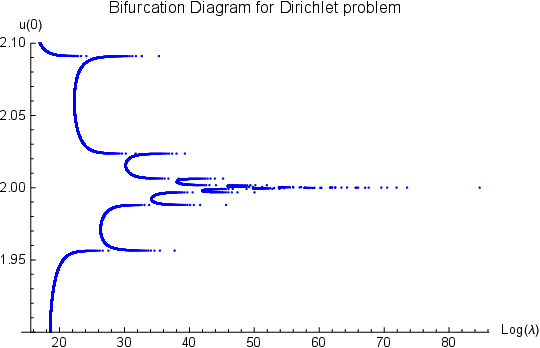}}
\caption{Infinitely many Dirichlet curves and infinitely many ground state solutions on both sides of $u(0)=2$}
\label{fig:3}
\end{figure}

The solutions of the initial value problem (\ref{d2}) with $u(0)=u_i$ and $u(0)=v_i$ provide a double infinite sequence of ground state solutions, see Figure \ref{fig:3}. The Lin-Ni solution in this case is $u(r)=\frac{2}{\sqrt[3]{36 r^2+1}}$, with $u(0)=2$. Thus we have two infinite sequences of ground states converging to the Lin-Ni solution, consistent with the result of I. Flores \cite{F}.
\medskip

\begin{figure}[h!]
\scalebox{1.2}{\includegraphics{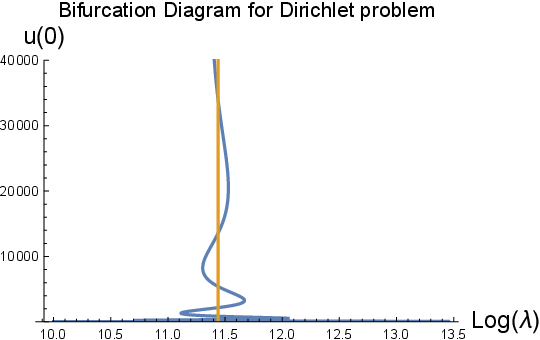}}
\caption{The curve approaching infinity}
\label{fig:2}
\end{figure}

For $u(0)$ large the solutions of (\ref{d1}) lie on a unique solution curve, as follows by Proposition \ref{prop:2}. This curve approaches infinity as $\la \ra \la _{\infty} \approx 92950$, making infinitely many turns around the line $\la =\la _{\infty}$, see Figure \ref{fig:2}.

\medskip

\noindent
The Figures \ref{fig:3} and \ref{fig:2} above were drawn in case $n=3$ and $p=4$.  Since $\la $ has to be very large for the features of the curve to come out, we draw $  u(0)$ versus $\ln \la$ in both cases (which required changing $r^2$ to $2 \ln r$ in the AppendTo command). The Figure \ref{fig:2} in particular shows that no a priori estimates are possible for the equation (\ref{d1}).
\medskip

The extraordinary complexity of the solution set  in any arbitrarily small neighborhood of the  Lin-Ni solution makes it interesting to study the linearized problem of  (\ref{d1}) 
\[
w''+\frac{n-1}{r}w'+\left(pu^{p-1}+(2p-1)u^{2p-2} \right)w=0 \,, \s w'(0)=0 \,, \; w(1)=0 
\]
at the Lin-Ni solution. We took $n=3$, $p=4$, so that  the linearized equation at  the Lin-Ni solution $u(r)=\frac{2}{\sqrt[3]{36 r^2+1}}$ becomes
\[
w''+\frac{2}{r}w'+\frac{96 \left(12
   r^2+5\right)}{\left(36
   r^2+1\right)^2} w=0 \,.
\]
Since $\frac{96 \left(12
   r^2+5\right)}{\left(36
   r^2+1\right)^2} \sim \frac{8}{9r^2}$ for large $r$,  solutions of the linearized equation are approximated by those of Euler's equation
\beq
\lbl{d10}
z''+\frac{2}{r}z'+\frac{8}{9r^2} z=0 \,.
\eeq
The general solution of (\ref{d10}) is 
\[
z=\frac{c_1 \cos
   \left(\frac{1}{6} \sqrt{23}
   \log
   (r)\right)}{\sqrt{r}}+\frac
   {c_2 \sin \left(\frac{1}{6}
   \sqrt{23} \log
   (r)\right)}{\sqrt{r}} \,.
\]
Hence we expect $w(r)$ to make infinitely many oscillations as $r \ra \infty$.
\medskip

Our numerical computations are supported by the following typical result on the discreteness and boundness of ground states.

\begin{prop}
There is a discrete sequence of ground state solutions of Lin-Ni equation (here n=3)
\beq
\lbl{d10+}
u''+\frac{2}{r}u'+u^4+u^7=0, \s u'(0)=0.
\eeq
These ground states have the maximal values  $u(0)$ tending to $u(0)=2$ from both above and below, and they separate an infinite sequence of  solution curves of the Dirichlet problem
\beq
\lbl{d10a}
u''+\frac{2}{r}u'+\la \left( u^4+u^7 \right)=0, \s u'(0)=0, \s u(1)=0.
\eeq
Moreover, for $u(0)$ large, all solutions of (\ref{d10a}) lie on a unique solution curve. 
\end{prop}

\pf
The Lin-Ni solution $u(r)=\frac{2}{\sqrt[3]{36 r^2+1}}$ is a ground state with $u(0)=2$, and it is of slow decay. It was shown in I. Flores \cite{F} that  there are infinitely many fast decay ground states. It is not possible to have an interval containing $u(0)=2$ giving ground states, since they would have to oscillate around the slow decay  Lin-Ni solution. Hence, the fast decay  ground states form a discrete sequence, and there are Dirichlet solution curves with $u(0)$ lying between them.
\medskip

By continuity, ground states of (\ref{d10+}) with $u(0)$ large would have to be close to singular solutions of (\ref{d10a}), and thus have a finite root. This is not possible, and hence large ground states do not exist. Since  ground states separate curves of solutions of (\ref{d10a}), it follows that solutions of (\ref{d10a}) with large $u(0)$ lie on a unique solution curve.
\epf

Our numerical computations show that along the curve  where $u(0) \ra \infty$, one has $\la \ra \la _{\infty} \approx 92967$, and $\la$ oscillates infinitely many times around the line $\la =\la _{\infty}$. In case $q>p_{JL}$, the Joseph-Lundgren exponent defined above, our numerical computations indicate that there are only at most finitely many oscillations around the vertical asymptote $\la = \la _{\infty}$. This result was conjectured in Y. Miyamoto \cite{M}, p. 676, and it appears to be still open.
\medskip

We now mention a related problem
\beq
\lbl{+*}
u''+\frac{n-1}{r}u'+\la \left( 1+u \right)^q=0, \s u'(0)=0, \s u(1)=0.
\eeq 
The following theorem was proved by D.D. Joseph and T.S. Lundgren \cite{JL}, and a simpler proof of a considerably more detailed result was given in P. Korman \cite{K}.

\begin{thm}
Assume that $\frac{n+2}{n-2}<q<p_{JL}$. Then all positive solutions of (\ref{+*}) lie on a unique continuous solution curve, starting at $\la=0$, $u(0)=0$, with $u(0) \ra \infty$, as $\la \ra \la _{\infty}=B$ (with $B$ defined in (\ref{15}), and with $q$ in place of $p$). As $\la \ra \la _{\infty}$, we have $u(r) \ra r^{-\frac{1}{q-1}}-1$, the singular solution.
\end{thm}

Let us briefly recall the proof of this theorem in \cite{K}. Using self-similarity of the problem (\ref{+*}), one can explicitly express the global solution curve of  (\ref{+*}) in the form
\[
\left(\la ,u(0) \right) = \left(t^2w^{q-1}(t),-1+\frac{1}{w(t)}  \right), \s 0 \leq t<\infty,
\]
where $w(t)$ is the solution of 
\[
w''+\frac{n-1}{r}w'+w^q=0, \s w(0)=1, \s w'(0)=0,
\]
which is well known to exist for all $t>0$, with $\lim _{t \ra \infty}w(t)=0$. Then it follows by the S. Chandrasekhar's Lemma \ref{lma:50} that $ \la _{\infty}=B$.
\medskip

Now suppose that $n=3$ and $q=7$, so that the nonlinear terms in (\ref{d10a}) and (\ref{+*}) are asymptotically the same as $u \ra \infty$. Are the solution curves similar for large $u(0)$? Qualitatively yes, but for (\ref{+*}), one has $\la _{\infty}=B_{\left(n=3,p=4 \right)} \approx 0.22$, much smaller than $\la _{\infty} \approx 92967$ for (\ref{d10a}). This sharp change in $\la _{\infty}$ occurs because $f(u)=(1+u)^7$ is not small for $u$ small, and hence the singular solution of (\ref{+*}) (with $\la$ scaled out) develops the first root much faster than for $f(u)=u^p+u^q$. It is an interesting open problem to search for a formula giving $\la _{\infty}$ for problems (\ref{i1}) that include  (\ref{d10a}).

\section{Singular solution for the general case}
\setcounter{equation}{0}
\setcounter{thm}{0}
\setcounter{lma}{0}

In this section we determine $\la _{\infty}$ at which the upper solution curve of
\beq
\lbl{s0}
u''+\frac{n-1}{r}u'+\la \left(u^p+u^q\right)=0, \s u'(0)=0, \s u(1)=0,  
\eeq
with 
\beq
\lbl{s0a}
1<p<\frac{n+2}{n-2}<q, \s n \geq 3,
\eeq
tends to infinity, and prove the existence of a singular solution at $\la=\la _{\infty}$. Singular solutions appeared first in the classical paper of H. Brezis and L. Nirenberg \cite{BN83}.
\medskip

In a number of special cases it is known that at some $\la=\la _{\infty}$ there exists a singular solution  of the equation in (\ref{s0}), which satisfies
\beqa
\lbl{s0b}
& u''+\frac{n-1}{r}u'+\la _{\infty} \left(u^p+u^q\right)=0, \\
& \lim _{r \ra 0} u(r)=\infty, \s u(1)=0, \nonumber
\eeqa
and solutions of  (\ref{s0}) tend to the singular solution as $\la \ra \la _{\infty}$, see C. Budd and J. Norbury \cite{BN}, F. Merle and L.A. Peletier \cite{me}, Y. Miyamoto \cite{M}, J. Dolbeault and I. Flores \cite{D}, see also P. Korman \cite{K4}. Here we present a simple proof of a more detailed result which is moreover suitable for numerical computations, and more flexible toward generalizations. Further, it was suggested in Y. Miyamoto \cite{M} that the Morse index of  the singular solution ``encodes" the oscillation properties of the solution curve of (\ref{s0}). In particular, if this index is infinite, the solution curve makes infinitely many turns as $u(0) \ra \infty$.
\medskip

We make a  change of variables
\beq
\lbl{snew}
 r = \frac{1}{\sqrt{\la}}t,
\eeq
and then switch back to the variable $r$, transforming this equation into
\beq
\lbl{s1}
u''+\frac{n-1}{r}u'+u^p+u^q=0\,.
\eeq
It follows that if $u(r)$ is the singular solution of (\ref{s1}), i.e. $\lim _{r \ra 0} u(r)=\infty$, and $\xi$ is its first root, then $\la _{\infty}=\xi ^2$.
Searching for the singular solution of (\ref{s1}), 
let
$
u=r^{\al}v $.
Substitute this into (\ref{s1}), then divide by $r^{\al-2}$ to get
\beq
\lbl{s2}
\s\s\s\s r^2v''+\al _1 rv'+\al _2v+r^{p_1}v^p+r^{q_1}v^q=0 \,,
\eeq
with $\al _1=2 \al+n-1$, $\al _2=\al (\al+n-2)$, $p_1=\al p-\al+2$, $q_1=\al q-\al+2$.
Choose $\al$: 
\beq
\lbl{s3}
\al =-\frac{2}{q-1} \,,
\eeq
the known growth rate of the singular solution, which also makes the power $q_1=0$. Obtain
\beq
\lbl{s4}
\s\s\s\s r^2v''+\al _1 rv'+\al _2v+r^{\gamma}v^p+v^q=0 \,,
\eeq
where $\al _1=-\frac{4}{q-1}+n-1$, $\al_2=-\frac{2}{q-1}\left(-\frac{2}{q-1}+n-2  \right)$. 
\beq
\lbl{s4a}
\gamma=2 \frac{q-p}{q-1}.
\eeq 
It is easy to check that under the condition (\ref{s0a})
\[
\al _1>0, \s \al _2<0,
\]

\noindent
{\bf Case 1.} Let us begin with the Lin-Ni equation, where $q=2p-1$, which turns out to be special. Then $\gamma=1$, and (\ref{s4}) takes the form
\beq
\lbl{s5}
\s\s\s\s r^2v''+\al _1 rv'+\al _2v+rv^p+v^{2p-1}=0 \,.
\eeq
This equation is singular at $r=0$ (although its linear part is similar to the Laplacian). It turns out that the value of $v''(0)$ is unbounded, except for a unique pair $(v(0),v'(0))$ of initial conditions at $r=0$, which we shall determine next. We are   searching for solution of  (\ref{s5}) in the form 
\beq
\lbl{s6}
v(r)=\sum _{n=0}^{\infty} a_n r^n=a_0+a_1r+a_2r^2+\cdots.
\eeq
Substitute this series into (\ref{s5}). In order for the constant term to vanish in (\ref{s4}), one needs
\[
a_0 \al _2+a_0^{2p-1}=0 \,,
\]
or $a_0={\left(-\al _2 \right)}^{\frac{1}{2p-2}}$. Equating to zero the coefficient in $r$ gives
\[
a _1 \al _1+a_1 \al _2+a_0^p+(2p-1)a_0^{2p-2}a_1 =0 \,,
\]
or $a_1=-\frac{a_0^p}{\al _1+\al _2+(2p-1)a_0^{2p-2}}=-\frac{a_0^p}{\al _1-2(p-1)\al _2}$. 

\medskip
Equating to zero the coefficient in $r^k$, $k \geq 2$, one can recursively calculate other $a_k$ through $a_0,a_1, \ldots,a_{k-1}$ as in \cite{K4}, and then prove the convergence of the series in (\ref{s6}) similarly to \cite{K4}, justifying the existence of a unique singular solution of (\ref{s2}) in the form $u(r)=r^{-\frac{2}{q-1}}v(r)$ with {\em analytic} $v(r)$. (Alternatively, the coefficients of the series in (\ref{s6}) can be calculated by repeated differentiation of the equation (\ref{s5}),  similarly to \cite{K4}.)
\medskip

Alternatively, for numerical computations we  used  $v(0)=a_0$ and $v'(0)=a_1$ as the unique set of  initial conditions for (\ref{s4}), for which one can use the NDSolve command in {\em Mathematica} to calculate the solution $v(r)$, and its first root $\xi$ (which is the same as the first root of the singular solution $u(r)$), so that $\la _{\infty}=\xi ^2$. To see that the root $\xi$ exists, observe that by continuity the singular solution of (\ref{s1}) is close to the solutions $u(r,\al)$ of the same equation with $u(0)=\al$, $u'(0)=0$ and $\al $ large. By Proposition \ref{prop:2a} we have positive bounds from below and above on the first root of $u(r,\al)$.
\medskip

In view of (\ref{snew}), we calculated the singular solution of (\ref{s0b}) as $t^{-\frac{2}{q-1}}v(t)$. Going back to $r$, by (\ref{snew}) we conclude that the singular solution is
\beq
\lbl{snew1}
u(r)=\la _{\infty}^{-\frac{1}{q-1}} r^{-\frac{2}{q-1}} v(\sqrt{\la _{\infty}}r).
\eeq

\noindent
{\bf \large Example 1 } For $n=3$, $p=4$ and $q=2p-1=7$, the form of the singular solution by (\ref{snew1}) is $u(r)=\la _{\infty}^{-\frac{1}{6}} r^{-\frac{1}{3}}v(\sqrt{\la _{\infty}}r)$, and a calculation shows that $v(0)=a_0=\frac{2^{1/6}}{3^{1/3}}$, $v'(0)=a_1=-\frac{1}{4 \, 6^{1/3}}$, and solving (\ref{s4}) (using the NDSolve command, with $h=$\$MachineEpsilon, $v(h)=a_0$, $v'(h)=a_1$) obtain $\ln \la_{\infty} \approx 11.44$, consistent with our numerical calculations of the largest solution curve $(\la,u(0))$ of (\ref{s0}), along which we had $u(0) \ra \infty$ and $\la \ra \la _{\infty}$. Observe that $\la _{\infty}$ is very large.
\medskip

If $p$ is other than the Lin-Ni exponent, it is natural to search for solution of (\ref{s4}) in the form (with $\gamma$ as defined by (\ref{s4a}))
\beq
\lbl{s8}
v(r)=\sum _{n=0}^{\infty} a_n r^{n \gamma}=a_0+a_1r^{ \gamma}+a_2 r^{2 \gamma}+ \cdots,
\eeq
which was suggested in C. Budd and J.  Norbury \cite{BN} (the equation (\ref{s10}) below is similar to (\ref{s5}), which motivates this form of solution). The Lin-Ni exponent $p=\frac{q+1}{2}$ turns out to be critical.
\medskip

\noindent
{\bf Case 2.} We consider next the case when $p$ is below the Lin-Ni exponent, $1<p<\frac{q+1}{2}$, $q>\frac{n+2}{n-2}$. Here $1<\gamma <2$.
Substitute this series into (\ref{s4}). In order for the constant term to vanish in (\ref{s5}), one needs
\[
a_0 \al _2+a_0^{q}=0 \,,
\]
or $a_0={\left(-\al _2 \right)}^{\frac{1}{q-1}}$. One could compute other $a_n$ in (\ref{s8}) similarly to Case 1, however since $\gamma>1$, we have two initial conditions for  solving (\ref{s4}), $v(0)=a_0$ and $v'(0)=0$, and we conclude the existence and uniqueness of solution of (\ref{s4}) as above, and use the NDSolve command for numerical computation. If $\xi$ is the first root of $v(r)$,  then as above $\la _{\infty}=\xi ^2$. Our numerical computations confirm the accuracy of this method. This provides a computational method. The form of solution in this case is given below.
\medskip

\noindent
{\bf Case 3.} Finally, assume that $\frac{q+1}{2}<p<\frac{n+2}{n-2}$ and $q>\frac{n+2}{n-2}$. Then $\gamma <1$ in (\ref{s8}), so that $v'(0)$ is not defined and the above approach does not work.
Making a change of variables $z=r^{\gamma}$ in (\ref{s4}), obtain (here $v=v(z)$)
\beq
\lbl{s10}
\gamma ^2 z^2 v''+\al _3 zv' +\al _2 v +zv^p+v^q=0,
\eeq
where 
\[
\al _3=\gamma (\gamma-1)+\al _1 \gamma.
\]
This equation is essentially the same as the  equation (\ref{s5}) for the Lin-Ni case, with different coefficients. We search for the solution of (\ref{s10}) in the form
\beq
\lbl{s12}
v(z)=a_0+a_1z+a_2z^2+\cdots .
\eeq
which is consistent with (\ref{s8}). Substitute this series into (\ref{s10}). In order for the constant term to vanish in (\ref{s10}), one needs
\[
a_0 \al _2+a_0^q=0 \,,
\]
or $a_0={\left(-\al _2 \right)}^{\frac{1}{q-1}}$. Equating to zero the coefficient in $z$ gives
\[
a _1 \al _3+a_1 \al _2+a_0^p+qa_0^{q-1}a_1 =0 \,,
\]
or $a_1=-\frac{a_0^p}{\al _3+\al _2+qa_0^{q-1}}=-\frac{a_0^p}{\al _3+(1-q)\al _2}$. We can now solve  (\ref{s10}) by either using the NDSolve command with $v(0)=a_0$ and $v'(0)=a_1$, or by computing $a_2,a_3,\ldots$ in (\ref{s12}). Then $v(r^{\gamma})$ is the solution of (\ref{s8}). If $z=\eta$ is the first root of $v(z)$, then $\xi =\eta ^{\frac{1}{\gamma}}$ is the first root for the solution of (\ref{s8}), and $\la _{\infty}=\xi ^2$.
\medskip

The same approach applies for Case 2 as well (where $\gamma >1$). We conclude the following result.

\begin{thm}
If $p$ and $q$ satisfy (\ref{s0a}), then the problem (\ref{s0b}) has a unique singular solution $(\la _{\infty},u(r))$, which is of the form $u(r)=r^{-\frac{2}{q-1}}v(r^{\gamma})$, where the function $v(z)$ is analytic, and $\gamma$ is defined by (\ref{s4a}). If $1<p<\frac{q+1}{2}$ and $q>\frac{n+2}{n-2}$ then  $\gamma >1$, and $v(r)$ is differentiable  at $r=0$.
If $p=\frac{q+1}{2}$ and  $q>\frac{n+2}{n-2}$ then  $\gamma =1$, and $v(r)$ is analytic  at $r=0$. If $\frac{q+1}{2}<p<\frac{n+2}{n-2}$ and $q>\frac{n+2}{n-2}$ then  $\gamma <1$, and $v(r)$ is only continuous (not differentiable) at $r=0$.
\end{thm}

In all three cases $v(z)$ is given by a convergent series near $z=0$. However its radius of convergence can be small, and it is typically less than $1$ according to \cite{K4}. Therefore for numerical calculations it is preferable to determine the initial conditions leading to $v(z)$, and calculate $v(z)$ by a computer algebra system (using the NDSolve command in {\em Mathematica}).
\medskip

The existence of singular solution was proved previously by Y. Miyamoto \cite{M} by a very complicated method (based on F. Merle and  L.A. Peletier \cite{me}), which did not produce the detailed information of the theorem above (this theorem also revealed the special role of the Lin-Ni equation).
\medskip

Our construction of singular solution leads to the following interesting conclusion.

\begin{prop}
Assuming that condition (\ref{s0a}) holds, let $W_0(r)$ be the singular solution of (\ref{s1}), while $w_0(r)$ is the singular solution of
\[
u''+\frac{n-1}{r}u'+u^q=0.
\]
Then $\lim _{r \ra 0} \frac{W_0(r)}{w_0(r)}=1$.
\end{prop}

\pf
By Lemma \ref{lma:lane}, $w_0(r)=B^{\frac{1}{q-1}}r^{-\frac{2}{q-1}}$. When constructing the singular solution, in all three case we had $W_0(r) \sim a_0r^{-\frac{2}{q-1}}$ as $r \ra 0$, with $a_0=(-\al _2)^{\frac{1}{q-1}}r^{-\frac{2}{q-1}}$. Since $-\al _2=B$, the proof follows.
\epf
\medskip

Finally, we discuss the role of singular solution for the global solution curve  of (\ref{s0}). Based on our computations we conjectured the following result.

\begin{prop} \lbl{prop:7}
Consider the problem (\ref{s0}) with $p$ and $q$ satisfying (\ref{s0a}). Assume that the unique singular solution of (\ref{s0}) occurs at $\la =\la _{\infty}$. Then $\la  \ra \la _{\infty}$, as $u(0) \ra \infty$ along the highest solution curve of (\ref{s0}), while the solution $u(r)$ tends to the singular solution of (\ref{s0}) at $\la =\la _{\infty}$.
\end{prop}

This result was proved in a number of special cases in \cite{BN}, \cite{me}, \cite{K}. The proof for the general case will appear in \cite{Kn}.
\medskip

\noindent
{\bf Acknowledgments}. It is a pleasure to thank Junping Shi and Yi Li for useful discussions.

\end{document}